\documentclass{amsart}[11pt,bezier,a4paper]

\usepackage[all]{xy}
\usepackage[dvips]{epsfig}
\usepackage{enumerate}
\usepackage[dvips]{graphicx}
\usepackage{psfrag}
\usepackage{array}
\usepackage{floatflt}
\usepackage{color}              % Need the color package
\xyoption{dvips}
\xyoption{rotate}

\usepackage{amssymb}

\usepackage[english]{babel}

\newtheorem{theorem}{Theorem}[section]
\newtheorem{deftheorem}{Definition-Theorem}[section]
\newtheorem{lemma}[theorem]{Lemma}
\newtheorem{e-proposition}[theorem]{Proposition}

\newtheorem{definition}[theorem]{Definition}
\newtheorem{remark}[theorem]{Remark}%{\it Remark\/}

\newtheorem{fact}[theorem]{Fact}
\def\og{\leavevmode\raise.3ex\hbox{$\scriptscriptstyle\langle\!\langle$~}}
\def\fg{\leavevmode\raise.3ex\hbox{~$\!\scriptscriptstyle\,\rangle\!\rangle$}}

\renewcommand{\AA}{\mathcal{A}}
\newcommand{\CC}{\mathcal{C}}
\newcommand{\PP}{\mathcal{P}}
\newcommand{\MM}{\mathcal{M}(\mathcal{A})}
\newcommand{\GG}{\mathcal{G}}
\newcommand{\UU}{\mathcal{U}}

\begin{document}
\bibliographystyle{amsplain}

% place in the next line the header (rubrique) chosen for your article,
% if you know it (you can also have 2, format : Header1/Header2
\centerline{}
%\begin{frontmatter}

% Title, authors and addresses

% use the thanksref command within \title, \author or \address for footnotes;
% use the ead command for the email address,
% and the form \ead[url] for the home page:
% \title{Title\thanksref{label1}}
% \thanks[label1]{}
% \author{Name\thanksref{label2}}
% \ead{email address}
% \ead[url]{home page}
% \thanks[label2]{}
% \address{Address\thanksref{label3}}
% \thanks[label3]{}
\selectlanguage{english}
\title[On a theorem of P. Deligne]{Combinatorial remarks\\ on a classical theorem of Deligne}

% use optional labels to link authors explicitly to addresses:
% \author[label1,label2]{}
% \address[label1]{}
% \address[label2]{}
% The [label1] can be suppressed if there is only one address for all authors

\selectlanguage{english}
\author[Emanuele Delucchi]{Emanuele Delucchi}
%\ead{delucchi@mail.dm.unipi.it}
%\author[authorlabel2]{Author Name2}
%\ead{author.name2@email.address2}

\address{Dipartimento di Matematica, Universit\`a di Pisa, Largo
  Bruno Pontecorvo 5, 56127 Pisa, Italia}
\thanks{The author acknowledges support of the Swiss National Science Foundation.}
%\keywords{Hyperplane arrangements, Order of regions, Partially ordered sets, lattices, $K(\pi,1)$-problem.}

%\address[authorlabel2]{Address2}
\email{delucchi@mail.dm.unipi.it}
% If you know the dates of reception, and acceptation you can put them now;
%  idem the name of the person presenting the Note

\medskip
%\begin{center}
%{\small Received *****; accepted after revision +++++\\
%Presented by £££££}
%\end{center}
%\begin{document}
\maketitle
\begin{abstract}
\selectlanguage{english}
% Text of abstract in English
We examine Deligne's classical proof of the asphericity of simplicial
arrangements from the viewpoint of the combinatorics of the poset of
regions of the arrangement. This turns out to be very
natural. In particular, we show that an arrangement is simplicial only
if it satisfies Deligne's property on positive paths, thus answering a
question posed by Paris in \cite{PaS}.
%{\it To cite this article: A.
%Name1, A. Name2, C. R. Acad. Sci. Paris, Ser. I 340 (2005).}

%\vskip 0.5\baselineskip

%\selectlanguage{francais}
% Text of abstract in French
%\noindent{\bf R\'esum\'e} \vskip 0.5\baselineskip \noindent
%{\bf Deux observations sur un th\'eor\`eme classique de Pierre Deligne. }
%Deligne a prouv\'e que le compl\'ementaire de la complexificaison de chaque
%arrangement simplicial est asph\'erique. On examine sa d\'emonstration
%%%%%du point de vue de la combinatoire de l'ordre partiel des chambres de
%5l'arrangement. On trouve que c'est une facon naturelle d'approcher le
%sujet. en particulier
%{\it Pour citer cet article~: A. Name1, A. Name2, C. R. Acad. Sci.
%Paris, Ser. I 340 (2005).}

\end{abstract}
%\end{frontmatter}

% now the Version franÁaise abrÈgÈe, if it exists
%\selectlanguage{francais}
%\section*{Version fran\c{c}aise abr\'eg\'ee}
% Text of your Version franÁaise abrÈgÈe here.
% Note you do not need to repeat here equations that you use in the
% main text - for example 'voir (3)' is quite acceptable.

%\selectlanguage{english}
% main text

\section{Introduction}
\label{intro}

%Let $V$ be a d-dimensional complex vector space. 
An arrangement of hyperplanes is a finite set of affine or linear codimension 1 subspaces of $\mathbb{C}^d$.
The arrangement induces a stratification of the ambient space by its hyperplanes
and their intersections. The poset of strata ordered by reverse inclusion is customarily
perceived as the combinatorial data of the arrangement. 
A famous
open question in arrangement theory is the so-called $K(\pi, 1)$-problem. An arrangement
is said to be $K(\pi, 1)$ if its complement in $\mathbb{C}^d$ is aspherical.  It is an open question
whether being $K(\pi, 1)$ is a combinatorial property in general.\\

A real arrangement of hyperplanes is called {\em simplicial} if the
maximal regions of the stratification it induces on real space (its
{\em chambers}) are
cones over simplices.
In his seminal paper \cite{deligne}, Deligne proved that the
complexification of a simplicial arrangement is $K(\pi,1)$. Deligne's proof
consists essentially in two steps. Assuming simpliciality of the
arrangement, he first derives a technical
property of the category of directed paths on the arrangement graph
(called `property D' by Paris in \cite{PaS}),
and then uses this property to show contractibility of the universal
cover of the complement.\\

Edelman introduced in \cite{ed1} a partial ordering of the chambers of
a real arrangement as a geometric generalization of the weak order on
Coxeter groups. Since general arrangements are not symmetric, this ordering depends on a choice of a `base
chamber', and the orders associated to different base chambers can
have quite different properties. The order-theoretic properties of
these posets were studied (e.g. in \cite{ed2,ed3, ed4}), and formalized in
the general framework of oriented matroids (see \cite[Chapter
4]{BLSWZ}, \cite{ERW}). 

The weak order of Coxeter groups is an example
of a combinatorial Garside structure. The construction of the complex
associated in \cite{besgroup, CMW} to any Garside structure can be
generalized to complexified hyperplane arrangements and leads to the
construction of `Garside-type' combinatorial models for the covers of
complexified arrangements (see \cite[Chapter 6]{PhD}). These models
are tiled by copies of the order complexes of the posets of regions.

Bj\"orner, Edelman and Ziegler \cite{BEZ} studied the
structure of the orderings of the regions of some combinatorially
defined classes of real arrangements. In particular, they show that a
real arrangement is simplicial if and only if the ordering of its
regions with respect to any base chamber is a lattice.\\

 We adopt this combinatorial point of view on simplicial
arrangements in examining Deligne's proof.  This turns out to be a very
natural way of formulating the argument. We prove that an arrangement
is simplicial only if it satisfies property D, thus showing that the
first part of Deligne's theorem is indeed an equivalence. This answers
a question posed by Paris (see \cite[p. 168]{PaS}).

Moreover, we see that the language of posets allows a very compact
proof of the contractibility of Garside-type models for the
universal cover of the complement of a simplicial arrangement.\\

We will begin by laying down the
combinatorial framework and introducing the main tools for our work in
Section \ref{comb}.  We prove that Deligne's
first step is an equivalence in Section \ref{eq}. Section
\ref{cones} proves contractibility of the garside-type model
of the universal cover starting from the formulation of property D in
terms of posets. We close the paper with Section \ref{last}, a short appendix
containing some considerations on a possible weakening of the
simpliciality condition.

\section{Combinatorics of real arrangements}\label{comb}
\subsection{Basics} Let $\AA:=\{H_i\}_{i=1,\ldots,n}$ denote an arrangement of linear
hyperplanes in $\mathbb{R}^d$. The complement of $\AA$ in
$\mathbb{R}^d$ is given by a set $\CC(\AA)$ of disjoint contractible components that we call {\em
  chambers} of $\AA$. We write $\CC(\AA)$ for the set of chambers of
$\AA$. Choose a point $p_C$ in the interior of every chamber
$C$. Given two chambers $C,C'\in\CC(\AA)$, we
say that an hyperplane $H\in\AA$ {\em separates} $C$ from $C'$ if the
segment joining $p_C$ with $p_{C'}$ intersects $H$. The set of
hyperplanes separating $C$ from $C'$ will be denoted by $S(C,C')$. Two
chambers $C,C'\in\CC(\AA)$ are said to be {\em adjacent} if there is
only one hyperplane separating them. The hyperplanes that separate a
chamber $C$ from its adjacent chambers are called {\em walls} of $C$. Since all hyperplanes are linear,
the arrangement is centrally symmetric with respect to the origin.Thus,
if $C\in\CC(\AA)$ then $-C\in\CC(\AA)$. Linearity of the hyperplanes
implies also that every chamber is a cone with the origin as apex. The base space of this cone is a convex $d$-polyhedron. If this polyhedron is a simplex
for every $C\in\CC(\AA)$, then $\AA$ is called {\em simplicial}.

The {\em complexification of $\AA$} is the arrangement
$\AA_\mathbb{C}$ obtained by considering the defining forms for the
$H_i$'s over $\mathbb{C}$. Let $\MM:=\mathbb{C`}^d\setminus\bigcup\AA_\mathbb{C}$ denote the complement
of $\AA_\mathbb{C}$.

\subsection{Partially ordered sets} In our considerations we will use some terminology and facts about
the combinatorics and topology of partially ordered sets (or, as we will say from now, {\em posets}) that we briefly recall. A more detailed
introduction can be found e.g. in \cite{Sta, bj}.

Let $P$ be a finite poset. Two elements $p,q\in P$ are called
comparable if either $p\geq q$ or $p<q$. 
Given $p\in P$ we define the subposets
$P_{\leq p}:=\{q\in P\mid q\leq p\}$, $P_{\geq p}:=\{q\in P \mid
q\geq p\}$, $P_{< p}:=P_{\leq p}\setminus\{p\}$, $P_{\geq
  p}:=P_{> p}\setminus \{p\}$. Any pair of comparable elements $p<q$
of $P$ determines an {\em interval} $[p,q]:=P_{\geq p}\cap P_{\leq
  q}$. We say that $q$ {\em covers} $p$, and write $q\gtrdot p$, if there is no element between
$p$ and $q$, i.e., if $[p,q]=\{p,q\}$.
A subset of $P$ consisting of pairwise comparable elements is called a {\em chain}. The
{\em length} of a chain is its cardinality. If every maximal chain of
$P$ has the same cardinality, the poset is called {\em
  graded} %. Equivalently, $P$
 and admits a {\em
rank function}, i.e., a function $\rho:P\rightarrow\mathbb{Z}$ that
is constant on minimal elements of $P$ and such
that $\rho(p)=\rho(q) +1$ whenever $p$ covers $q$. 
We say that $P$ is {\em bounded } if it
possesses a maximal and a minimal element (that are usually denoted by
$\hat{1}$ and $\hat{0}$, respectively). The poset $P$ is called a {\em
lattice} if for every pair of elements $p,q\in P$ the posets $P_{\geq
p}\cap P_{\geq q}$ and $P_{\leq
p}\cap P_{\leq q}$ are nonempty and bounded. In this case $p\wedge q
:= \min (P_{\geq
p}\cap P_{\geq q})$ is the unique minimal upper bound, called {\em join},
of $p$ and $q$. Similarly, $p\vee q:= \max (P_{\geq
p}\cap P_{\geq q})$ is the unique maximal lower bound, called {\em meet}, of $p$
and $q$.

The 
{\em order complex } of $P$, denoted by $\Delta(P)$, is the simplicial
complex given by the chains of $P$.  Note that we will make no explicit distinction between an abstract simplicial complex and its geometric realization.
If $P$ has a maximal element $\hat{1}$, the order
complex of $P$ is clearly a cone with apex $\{\hat{1}\}$ over the
space $\Delta(P_{<\hat{1}})$ and thus, in particular, contractible.

If $v$ is a vertex of a simplicial
complex $K$, the {\em star} of $v$ is the subcomplex given by all simplices
that contain $v$ and their boundaries. The {\em link} of $v$ is given
by all simplices of the star of $v$ that do not have $v$ as a
vertex.

\subsection{The arrangement graph} Let $\Gamma(\AA)$ denote the simple
graph on the vertex set $\CC(\AA)$ where two vertices are joined by an
edge if and only if the corresponding chambers are adjacent. The {\em
  arrangement graph} $G(\AA)$ is an oriented graph with the same set
of vertices (i.e., $\CC(\AA)$) and a pair of opposite oriented edges
between every two adjacent chambers.

The arrangement graph can be realized geometrically as the 1-skeleton
of the Salvetti complex, i.e., a CW complex that is homotopy equivalent to
$\MM$ (see \cite{sal1}). Thus, paths on $G(\AA)$ correspond naturally to topological
paths in $\MM$. The way the Salvetti complex is constructed implies
that any two directed paths of minimal length with the same beginning-
and endpoint are homotopic. Therefore we may write $(C\rightarrow C')$
for the equivalence class of the paths that are directed from $C$ to
$C'$ and have minimal length  (called {\em positive minimal paths}), and abuse terminology by referring to it
as to {\em the} positive minimal path from $C$ to $C'$. In fact, two
paths are homotopic in $\MM$ if and only if they are related by a
sequence of substitutions of equivalent positive minimal paths. 
Paths will be denoted by greek lowercase letters, and composed by concatenation.

The quotient of the free category on $G(\AA)$ with respect to the
relation generated by identifying any two positive minimal paths with
the same begin- and endpoint is the {\em category of positive paths}
$\GG^+(\AA)$. 
%This is a {\em germ} in the sense of Bessis
%\cite{besgroup}. 
Completion of $\GG^+(\AA)$ gives the {\em arrangement
groupoid} $\GG(\AA)$, which is clearly an instance of the fundamental
groupoid of $\MM$. For a precise account of this construction and
its significance for the modeling of the arrangement covers, see
\cite[Chapters 2, 4, 6]{PhD}.

\begin{remark}
{ The objects and facts of this section were already present and proved in the seminal work
by Pierre Deligne \cite{deligne}, where positive paths are called {\em
galeries}. We choose to adopt the above viewpoint because of the
convenience of the notation for our purposes.}
\end{remark}

\subsection{The order of regions} We now define a partial ordering of the set of
regions of a real hyperplane arrangement that was introduced by
Edelman \cite{ed1} (see also \cite{ed2,BEZ} for further study of this object).
%Fix a chamber $C_0$ of the real arrangement $\AA$, and call it the {\em base chamber} of our construction. The ordering we will define depends on the choice of $C_0$. 

\begin{definition}%\label{orderofregions} Let $\AA$ be a real arrangement of hyperplanes. For any two chambers $C_1,C_2\in \CC(\AA)$ let $$S(C_1,C_2):= \{H\in \AA \mid H\textrm{ separates } C_1 \textrm{ from }C_2\}.$$ \index{$S(C_1,C_2)$: the set of all hyperplanes that separate the chambers $C_1$ and $C_2$}
 let $\AA$ be a real arrangement of linear hyperplanes, and fix a base chamber $C_0\in \CC(\AA)$. We define the {\bf\em partial order $\PP_{C_0}(\AA)$ with base chamber $C_0$} on the set $\CC(\AA)$ by setting
$$ C_1 \leq_{C_0} C_2 \textrm{ if and only if } S(C_0, C_1) \subseteq S(C_0,C_2). $$  \end{definition}

One sees that the Hasse diagram of any $\PP_{C_0}(\AA)$ is given by
$\Gamma(\AA)$ after suitable choice of the ``bottom vertex''.
It is natural to ask about the order-theoretic properties of this
poset. For terminology and basic defititions on posets, see \cite{birkhoff}. 
First of all, from the above definition it is not hard to prove the following basic fact that we remark for later reference, pointing to \cite[Corollary 4.2.11]{BLSWZ} for a proof.
\begin{remark}\label{intervals}{ Let $\AA$ be a real linear arrangement, let $C, C_1, C_2 \in \CC(\AA)$ and suppose $C_1<_{C} C_2$. Then the interval $[C_1,C_2]\subset \PP_C$ is isomorphic to $\PP_{C_1}(\AA)_{\leq C_2}$. Thus, the structure of an interval is the same in all poset of regions where the interval is defined. }
\end{remark}
It is clear that, for any $C_0\in\CC(\AA)$, the poset
$\PP_{C_0}$ is bounded by $C_0$ and $-C_0$.
 Moreover, the cardinality of the sets $S(C_0,C)$ is a rank function for
  $\PP_{C_0}(\AA)$ by \cite[Proposition 1.1]{ed1}. In particular, the rank of
  $\PP_{C_0}$ equals the cardinality of $\AA$. Thus, the following Lemma
  gives a `local' sufficient condition for $\PP_{C_0}(\AA)$ to be a
  lattice.\vspace{2mm}

\begin{lemma}[Lemma 2.1 of \cite{BEZ}]\label{covering}
 Let $P$ be a bounded poset of finite rank such that, for any $p,q\in
 P$, if $p$ and $q$ both cover an element $w$ then the join $p\vee q$
 exists. Then $P$ is a lattice.\end{lemma}

This lemma is one of the ingredients of the proof of the following
characterization of simplicial arrangements in terms of their posets
of regions.

\begin{definition} Let $\AA$ be a real
  arrangement of linear hyperplanes and let $\CC(\AA)$ denote the set
  of its chambers. We say that $\mathbf{\AA}$ {\bf\em satisfies the strong lattice
  property} if $\PP_{C_0}$ is a lattice for every $C_0\in\CC(\AA)$.
\end{definition}
\begin{theorem}[Theorem 3.1 and 2.4 of \cite{BEZ}]\label{simpl} A real arrangement of linear hyperplanes $\AA$ is
 simplicial if and only if it satisfies the strong lattice property.
\end{theorem}

%The following lemma will allow us to get some control on the structure
%of the poset of regions when the base chamber is not simplicial.

\subsection{Topology of $\MM$}

%Recall from the introduction that in this paper we are mainly
%concerned with the asphericity of $\MM$. 
Since the complement of a
hyperplane arrangement in complex space is always connected, asphericity of $\MM$ is equivalent to contractibility of the
universal covering space. One possible way to approach the $K(\pi,1)$-problem is
therefore to construct combinatorially defined complexes that model
the homotopy type of the universal cover of $\MM$. This was indeed the
way taken by Deligne in \cite{deligne}: he considered a model for the
universal cover that was obtained by gluing together many copies of the
unit ball of $\mathbb{R}^d$ respecting the stratification given by
the arrangement. Later on, Paris made this point of view more explicit
and formulated Deligne's argument using a complex that, in the case of
a linear arrangement, lifted to the universal cover the simplicial
structure of the Salvetti complex (see \cite{PaS}). On the other hand,
if $\AA$ is the reflection arrangement of a finite irreducible Coxeter
group other complexes were studied, exploiting in different ways the
symmetry of this situation (see e.g. \cite{ChD,D,bestvina}). We want to
emphasize here the construction of Bestvina \cite{bestvina}, that was later
formulated in the more general context of Garside groups by Charney,
Meyer and Whittlesey \cite{CMW}, who described a universal cover
complex for Coxeter arrangements that is tiled by  order complexes of
the weak Bruhat order. This construction can be seen as a specially
symmetric case of the following complex, that models the homotopy type of any complexified arrangement of linear
hyperplanes and can be obtained by appropriately gluing copies of
the order complexes of all posets of regions associated to the
arrangement, as was shown in \cite{PhD}.\vspace{2mm}

\begin{deftheorem}[see Section 3.2 of \cite{PhD}] \label{defth}
Let $\AA$ be a complexified arrangement of linear hyperplanes, and fix
  a chamber $C_0\in\CC(\AA)$. We
  define a simplicial complex $\UU(\AA)$ which vertices are
  all morphisms of $\GG(\AA)$ that start at $C_0$ (i.e., all
  equivalence classes of paths on $G(\AA)$ that start at $C_0$). 
%on the vertex set
%  $V(\UU):=\Mor(\GG(\A))$. 
This simplicial complex is defined by declaring a
  set $\{\gamma_1, \dots ,\gamma_{d+1} \}$ of paths to be a simplex if and only if there are positive minimal paths
  $\alpha_1,\dots ,\alpha_d$ with 
 $\gamma_{i+1} = \gamma_{i}\alpha_i$ for all $i=1,\dots, d$ and
 $\alpha_1...\alpha_d$ positive minimal.\\
The complex $\UU(\AA)$ is homotopy equivalent to the universal cover
  of the complement $\MM(\AA)$ of the complexification of $\AA$.\end{deftheorem}

\subsection{Oriented matroids}

We point out that this section can be
phrased purely combinatorially in terms of the {\em oriented matroid}
of the arrangement (by saying {\em ``element''}  instead of
``hyperplane'' and {\em ``tope''} instead of ``region''). Thus, everything
can be defined for arbitrary arrangements of pseudospheres, though it is not clear what the
topological meaning of the constructions would be. For a comprehensive
introduction and a general reference to oriented matroids, see \cite{BLSWZ}.

\section{Necessity of the Strong Lattice Condition }\label{eq}

The first part in Deligne's proof of asphericity of simplicial
arrangements is devoted to show that, if the arrangement $\AA$ is simplicial, the morphisms of the positive
category $\GG^+(\AA)$ (i.e., the positive paths) can be written in a
particular normal form. Though it was recently referred to as 
`the Deligne normal form' (see e.g. \cite{bestvina, CMW}), we introduce it by rephrasing in our language Definition  of \cite{PaS}.\vspace{2mm}

\begin{definition} Let $\AA$ real arrangement of hyperplanes and fix
  $\bar{C}\in\CC(\AA)$. The arrangement satisfies {\em Property D} if for
  every positive path $\gamma$ starting at $\bar{C}$ there is a chamber
  $C_\gamma$ such that one can write $\gamma\sim \gamma'(C'\rightarrow
  C)$ for a positive path $\gamma'$ if and only if $C'< C_\gamma$ in
  $\PP_C$, where $C$ is the chamber in which $\gamma$ ends.
\end{definition}\vspace{2mm}

Paris asked in \cite{PaS} whether there are arrangements that satisfy
property D but are not simplicial. We will answer this question
negatively.\vspace{2mm}

First of all, we remark that, in view of Theorem \ref{simpl}, the result of the first part of Deligne's argument can be stated as follows.\vspace{2mm}

\begin{theorem}[Equivalent to Theorem 1.19 (iii) of \cite{deligne}]
If the arrangement $\AA$ satisfies the Strong Lattice Property, then it satisfies property D.
\end{theorem}\vspace{2mm}

Deligne's proof starts with the assumption of simpliciality. Our
remark is that, looking
at it with today's eyes, his argument has to spend quite a lot of work in deriving some
technical properties that are immediate consequences of the lattice
structure of the $\PP_{C}$s. In fact, the proof can be written entirely in
terms of posets of regions (see \cite{PhD}). From this combinatorial
point of view we can answer Paris' question as follows.
%can use the following well-known lemma that describes
%the structure of $\P_C(\AA)$ when it is not a lattice. 

\begin{theorem}
If the real arrangement $\AA$ satisfies property D then it satisfies
the Strong Lattice Condition and is therefore simplicial.
\end{theorem}

\noindent{\bf Proof:} We will argue by contraposition. Suppose that $\AA$ does
  not satisfy the Strong Lattice Condition, i.e., that there is a
  chamber $C_0$ such that $\PP_{C_0}$ is not a lattice. By Lemma
  \ref{covering}, this is only possible if  there are chambers $A$, $B$, $C$
  such  that $A$, $B$
  cover $C$ in the poset $\PP_{C_0}$ and the join $A\vee B$ does
  not exist in $\PP_{C_0}$. Since the interval
  $[C,-C_0]=(\PP_{C_0})_{\geq C}$ in $\PP_{C_0}$ is isomorphic to the
  interval $[C,-C_0]=(\PP_C)_{\leq -C_0}$ in $\PP_C$ (see Remark \ref{intervals}), we may from now
  on consider the situation in the latter poset, that is therefore
  also not a lattice. In particular, the chamber $C$ cannot be
  simplicial (Theorem 3.1 of \cite{BEZ}). Still, the following lemma
  tells us something about the structure of $\PP_C$ `near the bottom'.

\begin{lemma}[Lemma 4.4.4 of \cite{BLSWZ}, ``realizable version'']\label{bls444} Let $A$, $B$, $C$, $K$ be chambers of $\AA$, 
and suppose that $A$ and $B$ are atoms in the interval $[C,K]$ of $\PP_C(\AA)$. Then there exists a sequence of atoms $A=A_0, A_1, \ldots, A_k=B$ and a sequence of other elements $T_1, T_2,\ldots, T_k$ in $[C,K]$ such that $[C,T_i]$ is elementary and contains $A_i$ and $A_{i+1}$, for all $1\leq i \leq k$. If the chamber $C$ is simplicial, then $k=1$. (See Figure 1 (a))
\end{lemma}

%\begin{minipage}
\begin{figure}[h]
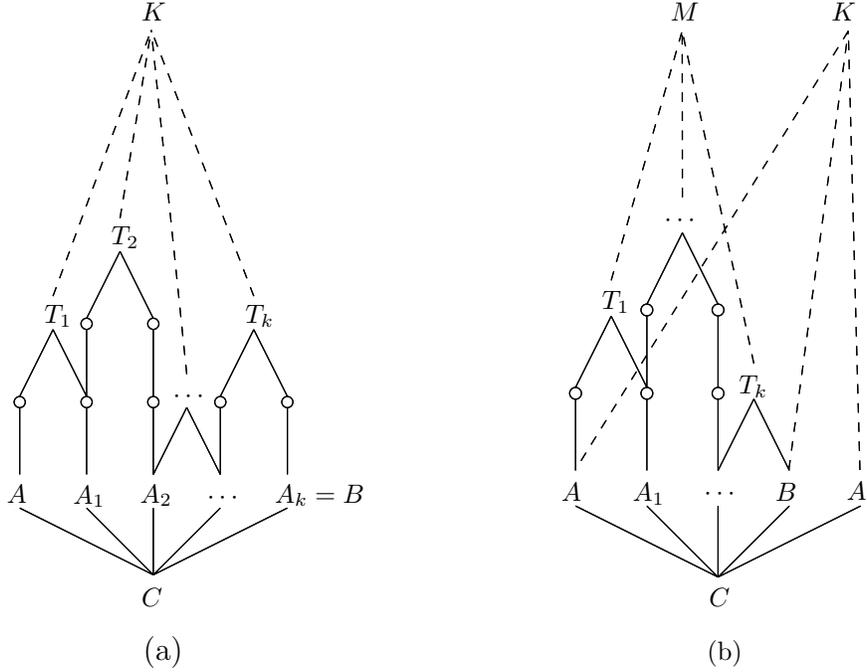
\label{lemma444}
\begin{center}
  \begin{picture}(0,0)%
    \includegraphics{fig_2_arr_3.pstex}%
  \end{picture}%
  \input{fig_2_arr_3.pstex_t}%
  
\caption{(a): The situation of lemma \ref{bls444}. %\newline
 (b): Our setup for the definition of $\gamma$.}
\end{center}
\end{figure}
%\end{minipage}
%\begin{minipage}
%
%\end{minipage}
In our setting, since the join of $A$ and $B$ does not exist, we have that, for any $K$, the minimal possible associated $k$ is at least $2$.  Let us choose an element $M$ among the minimal upper bounds of $A$ and $B$. We have the following situation (Figure 1 (b)): atoms $A=A_0, A_1,\ldots , A_k=B$ with the associated $T_1,T_2, \ldots T_k$, and all those elements are in the interval $[C,M]$.
%It has to be remarked that there might be other atoms $A'_j$ giving rise to the same situation as the $A_j$ (figure \thesection.3), with corresponding elements $T'_l$. Let $\mathscr{M}$ denote the set of all elements $M'$ that are minimal above all $T_l'$ for some sequence $(A'_j)_j$: without loss of generality we may choose a minimal element $M$ of $\mathscr{M}$, and suppose that our $A_j$ are the $A'_j$ giving rise to $M$. 
We will denote $H_i$ the unique element of $S(C,A_i)$.

Consider
$$\gamma:= (M\rightarrow C)(C\rightarrow A_1).$$

Clearly, $\gamma$ ends by both the positive paths $A\rightarrow C\rightarrow A_1$ and $B\rightarrow C\rightarrow A_1$. Let us show that {\em there is no chamber $K$ such that $\gamma$ ends with $K\rightarrow A_1$ and} $$S(A_1, K)\supset S(A,A_1)\cup S(B,A_1)=\{H_0, H_1,H_k\}.$$
Since lower intervals in lattices are closed under join, this shows
that $\AA$ does not satisfy property D and will therefore conclude the
proof.

Indeed, for such $K$ we would have $\{H_0, H_k\}\subset S(C,K)$ (because neither $H_0$ nor $H_k$ separate $A_1$ from $C$), therefore both $A$ and $B$ are atoms in the interval $[C,K]$. 
%An application of lemma \ref{bls444} to this situation yields a sequence of atoms$A'_0,\ldots, A'_l$ . 
So $K$ is clearly incomparable with $M$ in $\PP_C$: $K\not\geq M$ because $H_1\not\in S(C,K)$, and $K\not< M$ by minimality of $M$. This means that there is an hyperplane $H'\in S(M,K)\cap S(C,K)\setminus S(C,M)$, i.e., $H'$ separates $K$ and $C$ from $M$.

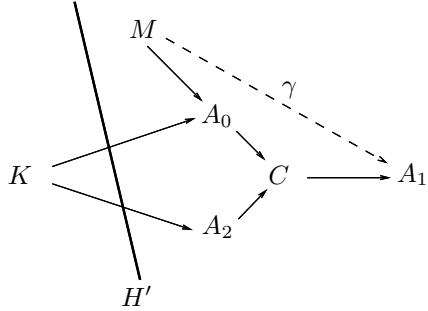
\begin{floatingfigure}[h]{7cm}\label{oursituation}
\begin{center}
  \begin{picture}(0,0)%
    \includegraphics{sketchpaths.pstex}%
  \end{picture}%
  \begin{picture}(0,0)%
\includegraphics{sketchpaths.pstex}%
\end{picture}%
\setlength{\unitlength}{3947sp}%
\begingroup\makeatletter\ifx\SetFigFont\undefined%
\gdef\SetFigFont#1#2#3#4#5{%
  \reset@font\fontsize{#1}{#2pt}%
  \fontfamily{#3}\fontseries{#4}\fontshape{#5}%
  \selectfont}%
\fi\endgroup%
\begin{picture}(2445,1973)(875,-1986)
\put(2513,-1179){\makebox(0,0)[lb]{\smash{\SetFigFont{10}{12.0}{\rmdefault}{\mddefault}{\updefault}{\color[rgb]{0,0,0}$C$}%
}}}
\put(1586,-1947){\makebox(0,0)[lb]{\smash{\SetFigFont{10}{12.0}{\rmdefault}{\mddefault}{\updefault}{\color[rgb]{0,0,0}$H'$}%
}}}
\put(1636,-265){\makebox(0,0)[lb]{\smash{\SetFigFont{10}{12.0}{\rmdefault}{\mddefault}{\updefault}{\color[rgb]{0,0,0}$M$}%
}}}
\put(3320,-1166){\makebox(0,0)[lb]{\smash{\SetFigFont{10}{12.0}{\rmdefault}{\mddefault}{\updefault}{\color[rgb]{0,0,0}$A_1$}%
}}}
\put(875,-1179){\makebox(0,0)[lb]{\smash{\SetFigFont{10}{12.0}{\rmdefault}{\mddefault}{\updefault}{\color[rgb]{0,0,0}$K$}%
}}}
\put(2094,-1496){\makebox(0,0)[lb]{\smash{\SetFigFont{10}{12.0}{\rmdefault}{\mddefault}{\updefault}{\color[rgb]{0,0,0}$A_2$}%
}}}
\put(2088,-798){\makebox(0,0)[lb]{\smash{\SetFigFont{10}{12.0}{\rmdefault}{\mddefault}{\updefault}{\color[rgb]{0,0,0}$A_0$}%
}}}
\put(2596,-607){\makebox(0,0)[lb]{\smash{\SetFigFont{10}{12.0}{\rmdefault}{\mddefault}{\updefault}{\color[rgb]{0,0,0}$\gamma$}%
}}}
\end{picture}

\end{center}
\end{floatingfigure}
Now, since $H_1$ is a wall of both $A_1$ and $C$, we have
 $(K\rightarrow A_1)=(K\rightarrow C)(C\rightarrow A_1)$. Therefore,
 if such a $K$ would exist, then $\gamma$ would consist of a positive
 path from $M$ to $K$ (thus crossing $H'$) followed by $(K\rightarrow
 A_1)$ (that crosses $H'$). 

But, by definition, $\gamma$ does not
 cross $H'$ since this hyperplane does not separate $M$ from $C$. 
%(see
% figure 2 for a sketch of the situation). 
This gives a contradiction: equivalent positive paths cross the same number of times every hyperplane (intuitively: $M\rightarrow K \rightarrow A_1$ `turns around' $H'$, while $\gamma$ does not).\hfill$\square$

\section{Contractibility of the universal cover}\label{cones}

Let us now turn to the second part of Deligne's proof of asphericity
of simplicial arrangements, where property D is used to show
contractibility of the universal cover. We want to phrase also this
step in combinatorial terms, using our complex $\UU(\AA)$ (recall Definition \ref{defth})

In analogy with Deligne's argument we have the following basic
observation, that is now standard.

\begin{lemma}[Proposition 6.4.4 of \cite{PhD}. See Proposition 2.14 of \cite{deligne}.]
Let $\AA$ be an arrangement of linear hyperplanes in
$\mathbb{R}^d$. Then $\UU(\AA)$ is contractible if the subcomplex
 $\UU^+$ given by the vertices that correspond to positive paths is contractile.
\end{lemma}

Now, to obtain the result it suffices to prove the following
statement. We will do this by using 

\begin{theorem}
Let $\AA$ be an arrangement of linear hyperplanes in
$\mathbb{R}^d$. If $\AA$ satisfies property D, then $\UU^+$ is contractible.
\end{theorem}
\noindent{\bf Proof: }  Let $\UU^+_m$ denote the subcomplex of $\UU^+$
given by the vertices that correspond to paths of edge-length at most
$m$. We will show that, for any $m>0$, $\UU^+_m$ retracts onto
$\UU^+_{m-1}$. 

Indeed, let $\gamma$ represent a vertex of $\UU^+_m\setminus
\UU^+_{m-1}$: it is a positive path of length $m$ that ends, say, in
the chamber $C$. Its link in $\UU^+_m$ is spanned by all vertices indexed
by positive paths $\gamma'$ such that
\begin{equation}\label{eq1}\gamma\sim\gamma'(C' \rightarrow C) \textrm{ and } C'\neq C,\end{equation} where
$C'$ denotes the chamber in which $\gamma'$ ends. 
Property D tells us that there is $C_\gamma\in\CC(\AA)$ such that
\ref{eq1} holds if and only if $C'\in I_\gamma:=[C,C_\gamma]_{>C}\subset\PP_{C}(\AA)$. 
%allowed chambers $C'$ are precisely the
%elements of $\PP_C$ that are below $C\gamma$ but strictly above
%$C$, i.e., the elements of the half-open interval $[C_\gamma,C]_{>C}$
It is easy to see that the subcomplex of $\UU^+$ spanned by the
elements of 
$I_\gamma$ is $\Delta(I_\gamma)$, which is contractible because $C_\gamma$ is a maximal element for
$I_\gamma$.  
%of $\PP_{C}$. By definition, the link of $\gamma$ is then the order complex
%$\Delta([C_\gamma,C]_{>C})$. Since the poset $[C_\gamma,C]_{>C}$ has a maximal element
%(i.e., $C_\gamma$), its order complex is a cone, thus contractible.

Thus, the link of $\gamma$ is a contractible subcomplex of
$\UU^+_{m-1}$, and we can then retract the star of $\gamma$ to it. Note that this process did not
involve any other vertex of $\UU^+_m\setminus\UU^+_{m-1}$ and can be
therefore be carried out successively for all vertices that correspond
to paths of length $m$.
Concatenation of the resulting retractions gives an explicit global retraction of $\UU^+_m$ onto $\UU^+_{m-1}$.\hfill$\square$

\section{Appendix: the Weak Lattice Property}\label{last}

Speaking about the order of regions as related to asphericity of
arrangements, it can not be omitted to mention a suggestive fact, obtained by
collecting results of \cite{BEZ,JP}.

\begin{fact}
Let $\AA$ be a linear arrangement of real hyperplanes. If $\AA$ is
simplicial, supersolvable or hyperfactored then there is
$C_0\in\CC(\AA)$ such that $\PP_{C_0}(\AA)$ is a lattice. 

We say that these classes of arrangements satisfy the `Weak Lattice Propery'.
\end{fact}

For background and definitions we refer to \cite{BEZ,JP,OT}. Here we
may only recall that simplicial and supersolvable arrangements are the
two combinatorially defined classes of arrangements that are up to now
known to be $K(\pi,1)$. Hyperfactored arrangements are a
generalization of supersolvable arrangements wor which the $K(\pi,1)$
property was never refuted.

It is clear that the Weak Lattice Propery does not imply asphericity
of arrangements (the arrangement $\AA_{-2}$ of Edelman and Reiner
\cite{ER} is not aspherical, but satisfies the Weak Lattice Property
by \cite[Theorem 3.2]{BEZ}). Nevertheless, it would be interesting to
investigate the significance of the structure of the poset of regions
for the asphericity of complexified arrangements. Some partial results in this sense can be found in
\cite{PhD}. We plan to expand on it in future work.

\bibliography{biblio}

\providecommand{\bysame}{\leavevmode\hbox to3em{\hrulefill}\thinspace}
\providecommand{\MR}{\relax\ifhmode\unskip\space\fi MR }
% \MRhref is called by the amsart/book/proc definition of \MR.
\providecommand{\MRhref}[2]{%
  \href{http://www.ams.org/mathscinet-getitem?mr=#1}{#2}
}
\providecommand{\href}[2]{#2}
\begin{thebibliography}{10}

\bibitem{besgroup}
David Bessis, \emph{Garside categories, periodic loops and cyclic sets}, ArXiv:
  \texttt{math.GR/0610778} (2006), 33 pp.

\bibitem{bestvina}
Mladen Bestvina, \emph{Non-positively curved aspects of {A}rtin groups of
  finite type}, Geom. Topol. \textbf{3} (1999), 269--302 (electronic).

\bibitem{birkhoff}
Garrett Birkhoff, \emph{Lattice theory}, third ed., American Mathematical
  Society Colloquium Publications, vol.~25, American Mathematical Society,
  Providence, R.I., 1979.

\bibitem{bj}
A.~Bj{\"o}rner, \emph{Topological methods}, Handbook of combinatorics, Vol.\
  1,\ 2, Elsevier, Amsterdam, 1995, pp.~1819--1872.

\bibitem{BEZ}
Anders Bj{\"o}rner, Paul~H. Edelman, and G{\"u}nter~M. Ziegler,
  \emph{Hyperplane arrangements with a lattice of regions}, Discrete Comput.
  Geom. \textbf{5} (1990), no.~3, 263--288.

\bibitem{BLSWZ}
Anders Bj{\"o}rner, Michel Las~Vergnas, Bernd Sturmfels, Neil White, and
  G{\"u}nter~M. Ziegler, \emph{Oriented matroids}, second ed., Encyclopedia of
  Mathematics and its Applications, vol.~46, Cambridge University Press,
  Cambridge, 1999.

\bibitem{CMW}
R.~Charney, J.~Meier, and K.~Whittlesey, \emph{Bestvina's normal form complex
  and the homology of {G}arside groups}, Geom. Dedicata \textbf{105} (2004),
  171--188.

\bibitem{ChD}
Ruth Charney and Michael~W. Davis, \emph{Finite {$K(\pi, 1)$}s for {A}rtin
  groups}, Prospects in topology (Princeton, NJ, 1994), Ann. of Math. Stud.,
  vol. 138, Princeton Univ. Press, Princeton, NJ, 1995, pp.~110--124.

\bibitem{D}
\bysame, \emph{The {$K(\pi,1)$}-problem for hyperplane complements associated
  to infinite reflection groups}, J. Amer. Math. Soc. \textbf{8} (1995), no.~3,
  597--627.

\bibitem{deligne}
Pierre Deligne, \emph{Les immeubles des groupes de tresses g\'en\'eralis\'es},
  Invent. Math. \textbf{17} (1972), 273--302.

\bibitem{PhD}
Emanuele Delucchi, \emph{Topology and combinatorics of arrangement covers and
  of nested set complexes}, Ph.D. thesis, ETH Zurich, 2006.

\bibitem{ERW}
P.~H. Edelman, V.~Reiner, and V.~Welker, \emph{Convex, acyclic, and free sets
  of an oriented matroid}, Discrete Comput. Geom. \textbf{27} (2002), no.~1,
  99--116.

\bibitem{ed4}
Paul~H. Edelman, \emph{The lattice of convex sets of an oriented matroid}, J.
  Combin. Theory Ser. B \textbf{33} (1982), no.~3, 239--244.

\bibitem{ed1}
\bysame, \emph{A partial order on the regions of {${\bf R}\sp{n}$} dissected by
  hyperplanes}, Trans. Amer. Math. Soc. \textbf{283} (1984), no.~2, 617--631.

\bibitem{ed3}
\bysame, \emph{Abstract convexity and meet-distributive lattices},
  Combinatorics and ordered sets (Arcata, Calif., 1985), Contemp. Math.,
  vol.~57, Amer. Math. Soc., Providence, RI, 1986, pp.~127--150.

\bibitem{ER}
Paul~H. Edelman and Victor Reiner, \emph{Not all free arrangements are
  {$K(\pi,1)$}}, Bull. Amer. Math. Soc. (N.S.) \textbf{32} (1995), no.~1,
  61--65.

\bibitem{ed2}
Paul~H. Edelman and James~W. Walker, \emph{The homotopy type of hyperplane
  posets}, Proc. Amer. Math. Soc. \textbf{94} (1985), no.~2, 221--225.

\bibitem{JP}
Michel Jambu and Luis Paris, \emph{Combinatorics of inductively factored
  arrangements}, European J. Combin. \textbf{16} (1995), no.~3, 267--292.

\bibitem{OT}
Peter Orlik and Hiroaki Terao, \emph{Arrangements of hyperplanes}, Grundlehren
  der Mathematischen Wissenschaften [Fundamental Principles of Mathematical
  Sciences], vol. 300, Springer-Verlag, Berlin, 1992.

\bibitem{PaS}
Luis Paris, \emph{Universal cover of {S}alvetti's complex and topology of
  simplicial arrangements of hyperplanes}, Trans. Amer. Math. Soc. \textbf{340}
  (1993), no.~1, 149--178.

\bibitem{sal1}
Mario Salvetti, \emph{Topology of the complement of real hyperplanes in {${\bf
  C}\sp N$}}, Invent. Math. \textbf{88} (1987), no.~3, 603--618.

\bibitem{Sta}
Richard~P. Stanley, \emph{Enumerative combinatorics. {V}ol. 1}, Cambridge
  Studies in Advanced Mathematics, vol.~49, Cambridge University Press,
  Cambridge, 1997.

\end{thebibliography}

% etc, etc

% The Appendices part is started with the command \appendix;
% appendix sections are then done as normal sections
% \appendix

% \section{}
% \label{}

% The Acknowledgements are an un-numbered section
%\section*{Acknowledgements}
% Acknowledgements text here

%\begin{thebibliography}{00}
% please try to use the bibitem system -
% the references should be in alphabetical order of authors' names.
% Articles with a single author first, author will 1 co-author next,
% then author with several co-authors;

% \bibitem{label}
% Texnnnnnnnnnnnnnnnt of bibliographic item
%\bibitem{bestvina}
%\bibitem{BEZ}
%\bibitem{BLSWZ}
%\bibitem{CMW}
%\bibitem{deligne}
%\bibitem{PaS}UC of SC
%\bibitem{PhD}
%\bibitem{ed1}
%\bibitem{ed2}
%\bibitem{sal1}
%\end{thebibliographynn}

\end{document}